\documentclass[11pt]{article}
\usepackage{amsmath,amssymb,amsfonts}
\usepackage{graphicx,subfigure}
\usepackage{graphics}

\usepackage{fullpage}
\newtheorem{prethm}{{\bf Theorem}}

\newenvironment{thm}{\begin{prethm}{\hspace{-0.5
em}{\bf.}}}{\end{prethm}}
\newtheorem{prepro}{{\bf Theorem}}

\newtheorem{precor}{{\bf Corollary}}

\newtheorem{cl}{Claim}

\newtheorem{preconj}{{\bf Conjecture}}

\newtheorem{preremark}{{\bf Remark}}

\newtheorem{prelem}{{\bf Lemma}}

\newtheorem{preproof}{{\bf Proof.}}

\title{Harmonious Coloring of Trees with Large Maximum Degree}

\author{{Saieed Akbari}\thanks{Department of Mathematical Sciences,
Sharif University of Technology, Tehran, Iran and School of Mathematics,
Institute for Research in Fundamental Sciences (IPM), P.O. Box 19395-5746,
Tehran, Iran. E-mail address: s$_-$akbari@sharif.ir. The research of this author was in part supported by a grant from IPM (No.90050214).
 } \and
{Jaehoon Kim} \thanks{Department of Mathematics, University of Illinois,
Urbana, IL, 61801, USA. E-mail address: kim805@illinois.edu.
Research of this author is partially supported by the
Arnold O. Beckman Research Award of the University of Illinois
at Urbana-Champaign.}
\and {Alexandr Kostochka}
\thanks{Department of Mathematics, University of Illinois,
Urbana, IL, 61801, USA and Sobolev Institute of Mathematics, Novosibirsk,
Russia. E-mail address: kostochk@math.uiuc.edu. Research of this author
is supported in part by NSF grant DMS-0965587, by 
the Ministry of education and science of the Russian Federation (Contract no. 14.740.11.0868) and by grant 09-01-00244-a
of the Russian Foundation for Basic Research.} }

\date{}

\begin{document}

\maketitle

\begin{abstract}
A  harmonious coloring of $G$ is a proper vertex coloring of $G$
such that every pair of colors appears on at most one pair of
adjacent vertices. The  harmonious chromatic number of $G$,
$h(G)$, is the minimum number of colors needed for a harmonious
coloring of $G$.  We show that if $T$ is a forest
of order $n$ with maximum degree $\Delta(T)\geq \frac{n+2}{3}$, then
$$h(T)=\left\{
  \begin{array}{ll}
    \Delta(T)+2, & \mbox{if $T$ has non-adjacent vertices of degree $\Delta(T)$;} \\
    \Delta(T)+1, & \mbox{otherwise.}
  \end{array}
\right.
$$
Moreover, the proof yields a polynomial-time
algorithm for an optimal harmonious coloring of such a forest.
\vspace{5mm}

\noindent {\small {\it Keywords}: Harmonious Coloring, Tree. }
{\small}

\vspace{-1mm}\noindent{\small {\it 2010 Mathematics Subject
Classification}}: {\small  05C05, 05C15.}

\end{abstract}

\vspace{1cm}

\section{ Introduction. }

Let $G$ be a simple graph. By $V(G)$ and $E(G)$ we denote the vertex set and the edge
set of $G$, respectively. A vertex of degree
$1$ in $G$ is called a {\it leaf}. A {\it harmonious
coloring} of $G$ is a proper vertex coloring of $G$ such that
every pair of colors appears on at most one pair of adjacent
vertices. The {\it harmonious chromatic number} of $G$, $h(G)$,
is the minimum number of colors needed for any harmonious
coloring of $G$. The first paper~\cite{fra} on harmonious coloring
appeared in 1982. However, the proper
definition of this notion is due to Hopcroft and Krishnamoorthy
\cite{hop}. 
Harmonious coloring of a graph is essentially an edge-injective homomorphism from a graph $G$ to a complete graph and the harmonious chromatic number of $G$ is the minimum order of a complete graph that admits such homomorphism from $G$. Paths and cycles are among the first graphs whose harmonious chromatic
numbers have been established \cite{fra}.  It was shown by
Hopcroft and Krishnamoorthy that the problem of determining the
harmonious chromatic number of a graph is NP-hard.
Moreover, Edwards and McDiarmid~\cite{kei1}
showed  that the problem remains hard even  restricted to
the class of trees. Since the problem is hard in the class of all trees, it makes sense to identify subclasses in which the problem is easier.

Since vertices at distance at most two in a graph $G$ must have distinct
colors in any harmonious coloring of $G$,  $h(G)\geq \Delta(G)+1$ for every graph $G$.
In \cite{afl}  it was shown that if $T$
is a tree of order $n$ and $\Delta(T)\geq \frac{n}{2}$, then
$h(T)=\Delta(T)+1$.  Moreover, the proof yields a polynomial-time
algorithm for an optimal harmonious coloring of such a tree.
We strengthen this result by finding a wider class of trees $T$ for which $h(T)=\Delta(T)+1$.

\begin{thm}\label{th1}
Let
\begin{equation}\label{l0}
\Delta \geq \frac{n+2}{3}.
\end{equation}
If $T$ is a forest of order $n$ with
$\Delta(T)=\Delta$, then
$$h(T)=\left\{
  \begin{array}{ll}
    \Delta+2, & \mbox{if $T$ has non-adjacent vertices of degree $\Delta$;} \\
    \Delta+1, & \mbox{otherwise.}
  \end{array}
\right.
$$
Moreover, there is a polynomial-time algorithm
 for an optimal harmonious coloring of such a forest.
 \end{thm}

In the next section we prove the lower bounds in the theorem  and show  that the bound $\Delta \geq \frac{n+2}{3}$ is sharp.
In the last two sections we prove the upper bounds in the theorem. 

Our notation is standard. In particular, for a graph $G$, $v\in V(G)$  and $W\subseteq V(G)$, $N_G(v)$ denotes the set of vertices
adjacent to $v$ in $G$, $d_G(v)=|N_G(v)|$, $N_G[v]=\{v\}\cup N_G(v)$, and $G[W]$ is the subgraph of $G$ induced by $W$.

\section{Lower bounds}

Since in each harmonious coloring $f$ of a graph $G$, the colors of all neighbors of a vertex $v$ are different 
and distinct from $f(v)$,
\begin{equation}\label{11}
\mbox{$h(G)\geq 1+\Delta(G)\qquad$ for every graph $G$.}
\end{equation}

\begin{cl}\label{cl1}
Let $k\geq 1$. If a graph $G$ contains two non-adjacent vertices, say 
$u_1$ and $u_2$,  of degree $k$, then $h(G) \geq k+2$.
\end{cl}
{\bf Proof.}
Suppose  that $G$ has a
harmonious $(k+1)$-coloring $f$ with colors in $A+\{\alpha_1,\ldots,\alpha_{k+1}\}$.
 Then by (\ref{11}), for each $j=1,\ldots,k+1$, the set $f^{-1}(\alpha_j)$
has a vertex in $N_G[u_1]$ and a vertex in $N_G[u_2]$.
 If $f(u_1)\neq f(u_2)$, then  the pair $\{f(u_1),f(u_2)\}$ appears on two
pairs of adjacent vertices: one pair in $N_G[u_1]$ and one pair in $N_G[u_2]$.
And if $f(u_1)= f(u_2)$, then 
for each $\alpha\in A-f(u_1)$, the pair $\{f(u_1),\alpha\}$ appears on two
pairs of adjacent vertices. So, $f$ is not harmonious. Thus $h(G)\geq k+2$.

\bigskip

Now for every $\Delta\geq 3$  we present a  tree $T_\Delta$ such that 
(i) $|V(T)|=3\Delta-1$, (ii) $\Delta(T)=\Delta$, (iii) $T$ has no non-adjacent vertices of degree $\Delta$, and (iv) $h(T)\geq \Delta+2$.
These examples show that the restriction (\ref{l0}) in Theorem \ref{l0} cannot be weakened.

Let $T_\Delta$ be obtained from a $4$-vertex path $(v_1,v_2,v_3,v_4)$ by adding $\Delta-1$ 
leaves adjacent to $v_1$, $\Delta-2$ 
leaves adjacent to $v_2$, and $\Delta-2$ leaves adjacent to $v_4$. Tree $T_4$ is depicted in Fig.~1.

\begin{figure}
\begin{center}\includegraphics[height=3cm]{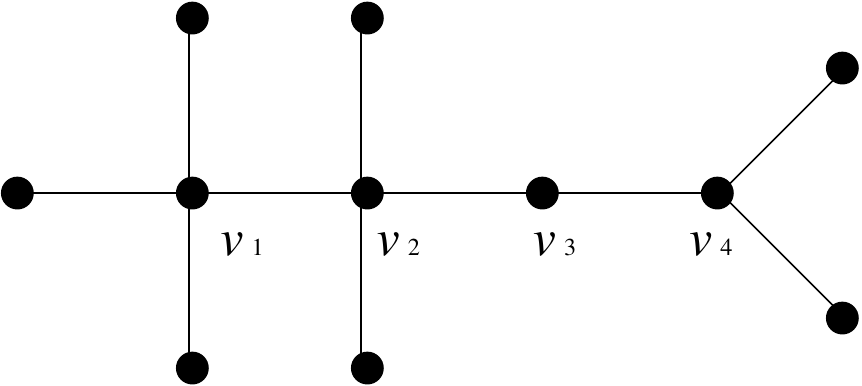}\end{center}
\caption{Tree $T_4$.}
\end{figure}

By construction, $T_\Delta$ is a tree with maximum degree $\Delta$ 
and $|V(T_\Delta)|=4+(\Delta-1)+(\Delta-2)+(\Delta-2) = 3\Delta-1 $. 
So, (i) and (ii) hold; and (iii) is also evident.
We establish (iv) by proving the following.

\begin{cl}\label{cl2}
$h(T_\Delta) =\Delta+2$.
\end{cl}
{\bf Proof.}
Suppose $f$ is a harmonious coloring of $T_\Delta$ with $\Delta+1$ colors.
We may assume that $f(v_1)=1, f(v_2)=2$, and $f(v_3)=3$.
Also we may assume that $f( N(v_1)-v_2 ) = \{3,4,\ldots, \Delta+1\}$.
Then $f(N(v_2) -v_1-v_3) = \{4,\ldots, \Delta+1\}$. Since $d_{T_\Delta}(v_1)=d_{T_\Delta}(v_2)=\Delta$,
no other vertices can be colored $1$ or $2$ in a harmonious $(\Delta+1)$-coloring.
 Thus, 
we may assume $f(v_4)=4$.  Then $f(N(v_4)-v_3) \subset \{5,6,\ldots,\Delta+1\}$ should hold. 
But $N(v_4)-v_3$ has $\Delta-2$ vertices and only $\Delta-3$ colors are available. 
Therefore we cannot complete the coloring with $\Delta+1$ colors. 
Thus $h(T_\Delta) > \Delta+1$.




\section{When $h(T)=\Delta+1$}
In this section, we present polynomial-time coloring procedures
 yielding that\\
 {\bf (*)} {\em if (\ref{l0}) holds and $T$ is an $n$-vertex forest with
$\Delta(T)=\Delta$ such that
 $T$ has no non-adjacent vertices of degree $\Delta$, then
$h(T)=\Delta+1$.}

First, observe that the statement holds for $\Delta\leq 2$: By  (\ref{l0}), $n\leq 3\Delta-2$.
So, if $\Delta\leq 1$, then $n\leq 1$, and so $h(T)\leq 1\leq 1+\Delta$. If $\Delta=2$, then
$n\leq 4$ and hence $T$ is a subgraph of the $4$-vertex path $P_4$ whose harmonious chromatic
number is $3=\Delta+1$. So, everywhere below 
\begin{equation}\label{au1}
\Delta\geq 3.
\end{equation}

Second, 
 let us check that it is enough to prove (*) for trees. Indeed,  if $T$ is a disconnected
$n$-vertex forest
satisfying (\ref{l0}) and (\ref{au1}) with no non-adjacent vertices of degree $\Delta$,
then by adding an edge connecting two leaves or isolated vertices from different components
of $T$, we again get a forest with these properties and fewer components. Thus in this section we will assume that $T$ is
a tree.

Let $v\in V(T)$ be a vertex of degree $\Delta$.
We will construct a harmonious coloring $f\,: V(T)\to V(K_{\Delta+1})$
step by step. The vertices of $H:=K_{\Delta+1}$ will by denoted by Greek
letters so that we do not mix them with the vertices of $T$.
We start from mapping $v$ and the $\Delta$ neighbors of $v$ in $T$ into
the all different $\Delta+1$ vertices of $H$.
Let $f(w)$ denote the color of $w$.
If $f(w)$ is not defined yet, then $w$ is  an {\em uncolored vertex}, otherwise  it is a
{\em colored vertex}.  

We consider several cases.

{\bf Case 0:}  $T$ consists of two stars with a path joining them.
 This case is straightforward.

{\bf Case 1:}  $v$ is the only vertex  of degree $\Delta$ in $T$, and
$T$ has no
vertices of degree $\Delta - 1$.
Suppose that we have already defined $f(w)$ for some vertices  $w\in V(T)$
(in particular, $f(w)$ is defined for $w\in N_T[v]$).
For $\alpha\in
V(H)$, let
$$d(\alpha):=\sum_{x\in
f^{-1}(\alpha)}d_T(x).$$
Also, we will say that vertices $\alpha$ and $\beta$ of $H$ are
{\em $T$-adjacent}, if there are $x\in f^{-1}(\alpha)$ and
$y\in f^{-1}(\beta)$ such that $xy\in E(T)$.
Our procedure will color one vertex at each step. It works as follows:

(a) Choose a vertex $w\in V(T)$ such that  $f(w)$ is defined
and $w$ has a neighbor $u$ for which $f$ is not defined and $u$ is not a leaf.
If there are no such vertex, then choose $u$ which is a leaf.

(b) If there is $\gamma\in V(H)-f(w)$ such that (i)  $\gamma$ is not $T$-adjacent to $f(w)$ and
(ii) $d(\gamma)+d_T(u)\leq \Delta$,
 then we let $f(u)$ be any $\gamma$ satisfying (i)--(ii) and
 go to (a) of the next step.

 (c) If no $\gamma\in V(H)-f(w)$ satisfies (i)--(ii), then we stop.

\medskip
We need to prove that we do not stop until we embed all $T$.
Note that after the initial coloring of $N_T[v]$,
 for every $\alpha\in V(H)$ we have $|f^{-1}(\alpha)|=1$, and hence
$d(\alpha)\leq\Delta$.

Suppose that we stopped in some step, before $f(x)$ was defined for
every $x\in V(T)$. This means that  at the moment of stopping, either every
$\gamma\in V(H)-f(w)$ is $T$-adjacent to $f(w)$ or
\begin{equation}\label{l1}
d(\gamma)+d_T(u)\geq \Delta+1\quad\mbox{for every
$\gamma\in V(H)-f(w)$ not $T$-adjacent to $f(w)$.}
\end{equation}
If the former holds, then since $f(u)$ is not defined yet,
at the moment of defining $f(w)$ we had already had
$d(f(w))+d_T(w)\geq \Delta+1$ and should have stopped then.
Thus some $\gamma\in V(H)-f(w)$ is not $T$-adjacent to $f(w)$,
and (\ref{l1}) holds. We may assume that $f(w)=\gamma_0$.
Let $\gamma_1,\ldots,\gamma_r$ be the
vertices of $H-\gamma_0$ not $T$-adjacent to $\gamma_0$, and
$\gamma_{r+1},\ldots,\gamma_{\Delta}$ be the
vertices of $H$ that are  $T$-adjacent to $\gamma_0$. And let $\gamma_{\Delta}=f(v)$. By the above,
$r\geq 1$.

By the choice of $u$, $u\neq v$. So in our case $d(u)\leq \Delta-2$. Thus by (\ref{l1}),
\begin{equation}\label{l2}
d(\gamma_i)\geq \Delta+1-d_T(u)\geq 3\quad\mbox{for every
$1\leq i\leq r$.}
\end{equation}

Since $f(v)$ is $T$-adjacent to every other vertex in $H$, according to our
rules, $f(x)\neq f(v)$ for every $x\neq v$.

\medskip

For every $W\subseteq V(T)$,
\begin{equation}\label{a1}
n-1=|E(T)|\geq \sum_{w\in W} d_T(w)-|E(T[W])|.
\end{equation}

We may assume that $d(\gamma_1)\leq d(\gamma_2)\leq \ldots \leq d(\gamma_r)$.
Let $$W:=\{v,u\}\cup f^{-1}(\{\gamma_0,\gamma_1,\ldots,\gamma_r\}).$$

{\em Case 1.1:} $r=1$.  Then, since $\gamma_0$ is $T$-adjacent to $\Delta-1$ vertices in $H$ and
$uw\in E(T)$, $d(\gamma_0)\geq \Delta$. By (\ref{l1}),
$d_T(u)+d(\gamma_1)\geq \Delta+1$. So since $\gamma_1$ is not $T$-adjacent to $\gamma_0$,
the graph  $T[W]$ has exactly $3$ edges, $uw, vx_0$ and $vx_1$, where $x_i$ is a neighbor of $v$ with $f(x_i)=\gamma_i$, for $i=0,1$. Thus using this $W$
  in  (\ref{a1}), we have
 $$|E(T)| \geq d_T(v)+d(\gamma_0)+d_T(u)+d(\gamma_1)-3\geq \Delta+\Delta+(\Delta+1)-3=3\Delta-2.$$
 So, $3\Delta-2 \leq n-1$, i.e.,  $\Delta \leq \frac{n+1}{3}$, a contradiction. \\

{\em Case 1.2:} $r=2$.  Similarly to Case 1.1,
 $d(\gamma_0)\geq \Delta-1$ and
$d_T(u)+d(\gamma_1)\geq \Delta+1$. Now $\gamma_1$ and $\gamma_2$ are not $T$-adjacent to $\gamma_0$.
So,   graph $T[W]$ has at most
 $5$  edges.  Thus, using this $W$ in  (\ref{a1}), we have (using also (\ref{l2}) to estimate $d(\gamma_2)$)
 $$|E(T)| \geq d_T(v)+d(\gamma_0)+d_T(u)+d(\gamma_1)+d(\gamma_2)-5\geq \Delta+(\Delta-1)+(\Delta+1)+3-5=3\Delta-2,$$
a contradiction as in Case 1.1.

{\em Case 1.3:}  $r\geq 3$. 
 Now $d(\gamma_0)\geq \Delta-r+1$ and
 $\gamma_1,\ldots,\gamma_r$ are not $T$-adjacent to $\gamma_0$. So,
 \begin{equation}\label{a2}
 |E(T[W])| \leq (1+r)+1+\sum_{i=1}^r\frac{d(\gamma_i)-1}{2}
 \end{equation}
   (here $r+1$ counts the edges incident with $v$, $1$ stands for the edge $uw$ and
 $\sum_{i=1}^r\frac{d(\gamma_i)-1}{2}$ estimates from above the number of edges both ends of which are in
 $\{\gamma_1,\ldots,\gamma_r\}$).
 So by  (\ref{a1}),
$$n-1= |E(T)| \geq
d_T(v)+d_T(u)+\sum_{i=0}^{r}d(\gamma_i)-r-2-\sum_{i=1}^{r}\frac{d(\gamma_i)-1}{2}$$
$$ \geq\Delta + (\Delta-r+1) + (d_T(u) + d(\gamma_1))-d(\gamma_1)-r-2 + \sum_{i=1}^{r} \frac{d(\gamma_i)+1}{2} $$
$$ \geq (2\Delta-r+1) + (\Delta+1)-d(\gamma_1)-r-2 +r\frac{d(\gamma_1)+1}{2}=
 3\Delta -2r + (r-2)\frac{d(\gamma_1)}{2}+\frac{r}{2}$$
 $$
 = 3\Delta -3 + (d(\gamma_1)-3) \frac{r-2}{2} \geq 3\Delta -3 + \frac{d(\gamma_1)-3}{2}.$$
 Thus if $d(\gamma_1)\geq 4$ or if  (\ref{a2}) is a strict inequality, or if
$d(\gamma_0)>\Delta-r+1$,
 then we have $n-1>3\Delta-3$,
which yields $\Delta\leq \frac{n+1}{3}$, a contradiction.
So, by (\ref{l2}) we may suppose that  $d(\gamma_1)=3$, $d(\gamma_0)=\Delta-r+1$, and  (\ref{a2}) holds with equality.
In particular,  by (\ref{l2}),  $d_T(u)=\Delta-2$.
 Since $r\geq 3$, we have $\Delta\geq 1+r\geq 4$, and so
 $u$ is not a leaf.
Since we do not have Case 0, there is a leaf $l$  not adjacent to $u$ and not adjacent to $v$.
 Thus,  according to our rule (a), $l$ is not colored yet. 
Since $d(\gamma_0)=\Delta-r+1$, and  (\ref{a2}) holds with equality,
$l$ is adjacent neither to any vertex in $f^{-1}(\gamma_0)$ nor
 to any vertex in $f^{-1}(\{\gamma_1,\ldots,\gamma_r\})$. Hence the right-hand side of (\ref{a1}) does not count the edge
incident with $l$. So, we have $n-2\geq 3\Delta -3$, a contradiction to $\Delta>\frac{n+1}{3}$.

Therefore we do not stop until we color all the vertices in $T$.

\medskip
Let $(u_1,\ldots,u_n)$ be an ordering of the vertices of $T$ such that $u_1=v$ and
$d_T(u_1)\geq d_T(u_2)\geq \ldots\geq d_T(u_n)$. In these terms, Case~1 was the case $d_T(u_2)\leq\Delta-2$.
Let $t$ be chosen so that $d_T(u_t)\geq 2$ and $d_T(u_{t+1})=1$.

{\bf Case 2:}  
$d_T(u_3)\geq\Delta-1$. Let $W':=\{u_1,\ldots,u_t\}$. Since $T$ is connected, $T[W']$ is also connected. Then
$$\sum_{w\in W'}d_T(w)\geq \Delta+(\Delta-1)+(\Delta-1)+2(t-3)\quad\mbox{and}\quad |E(T[W'])|= t-1.$$
So by  (\ref{a1}), 
$$n-1=|E(T)|\geq 3\Delta-2+2(t-3)-(t-1)\geq (n+2)+t-7=n+t-5.$$
It follows that $t\leq 4$, and that if $t=4$, then $d_T(u_2)=\Delta-1$ and $d_T(u_4)=2$.

{\em Case 2.1:} $t=3$. The only $3$-vertex tree is the 3-vertex path. So, 
$T[W']$ is the path $(w_1,w_2,w_3)$. By (*), we may assume that $w_3=u_3$.
For $i=1,2,3,$ we let $f(w_i)=\gamma_{i-1}$. We place the leaves adjacent to $w_1$ into
any $d_T(w_1)-1$ vertices in $V(H)-\gamma_0-\gamma_1$,  the leaves adjacent to
  $w_2$ into
any $d_T(w_2)-2$ vertices in $V(H)-\gamma_0-\gamma_1-\gamma_2$, and  the 
$\Delta-2$ leaves adjacent to $w_3$ into
the  vertices in $V(H)-\gamma_0-\gamma_1-\gamma_2$.

{\em Case 2.2:} $t=4$. As it was mentioned, in this case $d_T(u_2)=d_T(u_3)=\Delta-1$ and $d_T(u_4)=2$.
Since $T[W']$ is connected and $d_T(u_4)=2$, we may assume that $u_3$ is adjacent either to $u_1$ or
to $u_2$.  For $i=1,2,3,4,$ we let $f(u_i)=\gamma_{i-1}$. Then for $j=1,2,3,4,$ we place the 
leaves adjacent to $u_j$ into
the vertices of $H-\gamma_0-\ldots-\gamma_{j-1}$ not occupied by the neighbors of $u_j$ in $W'$. 
We can do it for $j=1,2$, since $d_T(u_j)= 1+\Delta-j$ for these $j$. And $u_3$ was chosen so that $\gamma_2=f(u_3)$
is $T$-adjacent to $\{\gamma_0,\gamma_1\}$. Finally, since $T[W']$ is connected, $u_4$ has at most one 
adjacent leaf.
So if $\Delta\geq 4$ or $u_4$ has no adjacent leaves, then we are done. Thus we need only to handle the situation
when $\Delta=3$ and each vertex of degree $2$ in $T$ has an adjacent leaf.
  Then $T[W']=K_{1,3}$. 
For $i=2,3,4$, let $w_i$ be the leaf adjacent to $u_i$.
In this case, we 
again let $f(u_i)=\gamma_{i-1}$ for $i=1,2,3,4$ and then let $f(w_2)=\gamma_2$, $f(w_3)=\gamma_3$,
and $f(w_4)=\gamma_1$.

\medskip
{\bf Case 3:}  $d_T(u_2)=\Delta$. Under Condition (*), $vu_2\in E(G)$, and
by Case 2, $d_T(u_3)\leq \Delta-2$.  
Since $vu_2\in E(G)$, $f(u_2)$ was defined at the first step.
Then  we can  apply the procedure of Case 1, and the argument goes through
since $d_T(u_3)\leq \Delta-2$.

\medskip
The only case, we have not yet considered is:

{\bf Case 4:} $d_T(u_2)= \Delta-1$ and $d_T(u_3)\leq \Delta-2$.
Let $P=(v_1, \ldots,v_q)$ be the path in $T$ connecting $v_1=v$
 with $v_q=u_2$.  Suppose $v$ has exactly $p$ adjacent non-leaves in $T$.
 We claim that 
 \begin{equation}\label{au3}
 q+p\leq \Delta+2,
\end{equation}
  since otherwise
 $$n\geq q+(p-1)+(d_T(v)-1)+(d_T(u_2)-1)\geq (\Delta+3-1)+(\Delta-1)+(\Delta-1-1)=3\Delta-1,$$
a contradiction to (\ref{l0}).
 
By (\ref{au3}), we can place all the vertices of $P$ and all remaining non-leaf neighbors of $v$ into
distinct vertices of $H$. After that, we place the leaves adjacent to
 $v$ into distinct vertices of $H$
not containing $v$ or its neighbors. Then we again apply the procedure  of Case 1
and the argument goes through
since $d_T(u_3)\leq \Delta-2$.

\section{Finishing the proof}
The only situation not covered in the previous section is that $T$ has
   non-adjacent vertices $v$ and $z$  of degree $\Delta$. 
We add a vertex $w$ to $T$ and make $w$ adjacent to $v$ to get a tree $T'$ with maximum degree $\Delta+1$. Then we may apply Case 1 to $T'$ and color $T'$ with $\Delta+2$ colors. This harmonious coloring of $T'$ gives a harmonious coloring of $T$ with $\Delta+2$ colors. This completes the proof of Theorem 1.

\bigskip
{\bf Remark.} Since we colored vertices one by one with no recolorings, and the choice
of every next vertex took polynomial time,  the total time taken by our algorithm is
polynomial.


\bigskip

{\bf Acknowledgments.} The authors would like to thank Bill Kinnersley for the helpful
comments and a referee for the idea of adding a vertex in Section 4.


\end{document}